\documentclass[12pt]{article}
\usepackage{epsfig}
\usepackage{amsmath}
\usepackage{amssymb}
\usepackage{color}
\normalsize \baselineskip 16pt
\newcommand {\beq} {\begin{equation}}
\newcommand {\eeq} {\end{equation}}

\begin{document}
\title{New Solutions for Slow Moving Kinks in a Forced Frenkel-Kontorova Chain}
\author{Phoebus Rosakis \\
\small Department of Applied Mathematics\\[-0.8ex]
\small University of Crete\\[-0.8ex]
\small Heraklion 71409, Greece\\[-0.8ex]
\small \texttt{rosakis@tem.uoc.gr}\\
\and
Anna Vainchtein\\
\small Department
of Mathematics\\[-0.8ex]
\small University of Pittsburgh\\[-0.8ex]
\small Pittsburgh, PA 15260, USA\\[-0.8ex]
\small \texttt{aav4@pitt.edu}\\
}

\maketitle

\begin{abstract} \noindent We construct new traveling wave solutions of moving kink type for a modified, driven, dynamic  Frenkel-Kontorova model, representing dislocation motion under stress.  Formal solutions known so far are inadmissible for velocities below a threshold value. The new solutions fill the gap left by this loss of admissibility. Analytical and numerical evidence is presented  for their existence; however, dynamic simulations suggest that they are probably unstable.
\end{abstract}

\section{Introduction}
The first study of dislocation motion through a discrete model was performed by Atkinson and Cabrera \cite{AC65}, followed by Celli and Flytzanis \cite{CF} and Ishioka \cite{IS71}.   Atkinson and Cabrera \cite{AC65} utilize a variant of the discrete sine-Gordon equation,  which arises in the dynamic version of the Frenkel-Kontorova  model.   They seek traveling wave solutions corresponding to uniformly moving dislocations under applied stress, represented by a constant forcing term.  To the present day, analytical  progress  on the forced problem has been limited, while some numerical conclusions have been drawn by Peyrard and Kruskal \cite{PK84}. In \cite{AC65}, the trigonometric potential of the Frenkel-Kontorova model, responsible for the nonlinearity in  the sine-Gordon equation, is replaced by a  piecewise smooth potential with quadratic wells.  A similar choice is made in  \cite{CF,IS71}.  The resulting equation for traveling waves reduces to  a linear one, provided the solution satisfies an  \emph{admissibility condition}. This requires that  to the right of a single transition point (the dislocation core), the solution take values entirely in one of  two quadratic valleys of the potential, and in the other valley to the left.  Solutions of the reduced linear equation are found semi-analytically.  Traveling waves are found to exist only if the stress and speed satisfy an algebraic relation, the \emph{kinetic relation} of the dislocation.  A remarkable prediction of this relation was that dislocations exceed the speed of shear waves at sufficiently high stress. This was recently confirmed experimentally \cite{NMR11}.  Another  feature is the presence of multiple singularities and discontinuities  in the graph of stress versus velocity, located at a sequence of  special resonance velocities that accumulate at zero; see  Fig.~\ref{fig:ACkinetics}. However, as emphasized by Earmme and Weiner  \cite{EW77}  (see also \cite{KT03}), below a \emph{threshold velocity} $V_0$, solutions of the linear problem violate the admissibility condition; this issue was not fully recognized in \cite{AC65}. This rules out solutions in the entire velocity interval from $0$ to  $V_0$  as inadmissible; this includes all singularities, since  all resonance velocities are below the threshold velocity.

The problem of existence of traveling waves at  velocities below the threshold value has remained open to the best of our knowledge; this also applies to the results of   \cite{CF}. This is the main issue addressed in the present paper.

We present compelling analytical and numerical evidence that  for velocities $V$ below the threshold value $V_0$, there is a new type of traveling wave solution $u_n(t)=u(\xi)$, $\xi=n-Vt$, of equation \eqref{eq:dyn} below, which actually \emph{violates} the usual admissibility conditions, but satisfies a suitable generalization. The model employs a piecewise smooth two-well potential $\Phi(u)$ with two quadratic branches meeting at  the \emph{spinodal value} $u=0$. The  usual admissibility conditions require that the solution  vanish (go through the spinodal value)  at precisely one point, say $\xi=0$, where it transitions between the two quadratic valleys of  $\Phi$, so that it must be strictly monotone in the neighborhood of the transition point. Instead, the new solutions we find
are equal to zero on an entire finite  interval, and lie in different potential valleys on either side of this interval.

The motivation for considering such solutions comes from a more elaborate model \cite{V10a} where $\Phi(u)$ is piecewise quadratic but continuously differentiable. Its graph consists of three parabolas, two convex ones separated by a concave one, defined on the \emph{spinodal range} of $u$ values.   Any kink solution $u(\xi)$  that transitions between the two convex potential valleys must therefore take values in the spinodal range for $\xi$ in some interval, say $[-z,z]$. Here $z>0$ is an unknown of the problem and is found to depend on the size of the spinodal range. Now $\Phi'(u)$ is trilinear with two increasing branches and a decreasing one between them. When the slope of the latter is treated as a parameter and approaches $-\infty$, so that the spinodal range tends to degenerate to the point $u=0$, a surprising observation is made in \cite{V10a}: at velocities below the threshold value $V_0$, the interval $-z\le\xi\le z$ where $u(\xi)$ takes values in the spinodal range does not shrink to a point, as one might expect. Rather, $z$ approaches a positive value in the limit.  The  limiting potential is the piecewise smooth biquadratic one considered in \cite{AC65,CF,IS71} where admissible  low-velocity kinks have not been found so far.  It is thus natural to  directly try  the new type of solutions mentioned in the previous paragraph, that vanish (equal the degenerate spinodal value) over an interval $[-z,z]$, where $z>0$, in the context of the Atkinson-Cabrera problem.

While we do not rigorously prove their existence, we are able to construct the new solutions semi-analytically; we present approximate analytical and also numerical evidence of their existence. We find that they bifurcate from  the classical Atkinson-Cabrera solutions, at precisely the threshold velocity $V_0$ below which the latter become inadmissible. The new solutions exist essentially over the entire  velocity range $0<V<V_0$; thus they seem to close the entire gap left by loss of admissibility of the Atkinson-Cabrera solutions.

The kinetic relation between the applied stress $\sigma$ (constant forcing term) and velocity $V$ of  moving kinks (dislocations) described by the new solutions  is entirely different from the (inadmissible)  one obtained in \cite{AC65} over  the velocity range $0<V<V_0$; see
Fig.~\ref{fig:new_vs_AC} where the two are compared. Unfortunately, its physical significance is in question, since the new solutions appear to be unstable. This is suggested by numerical simulations of the Frenkel-Kontorova chain dynamics, where kinks either stop or move with velocities above $V_0$.  Recently a similar trend was observed in experiments that measured kinetic relations for dislocations in actual two-dimensional plasma crystals \cite{NMR11}. Our results are consistent with the findings of \cite{V10a}, where slow kinks are apparently unstable when the  spinodal range is sufficiently narrow. In contrast, \emph{some} slow kinks do become stable if the spinodal range is sufficiently wide in \cite{V10a}, as is also observed in \cite{EW77,PK84}.

Our results clarify an important issue  associated with the dynamics of dislocation-like traveling wave solutions in a Frenkel-Kontorova chain. Based on a semi-analytical approach, they strongly suggest that there are such solutions for the whole range of velocities.  The role played by the two-well potential is also elucidated.  Use of the biquadratic non-smooth potential  with degenerate spinodal range (consisting of a point) yields low-velocity solutions that are apparently unstable.  In contrast, the results in \cite{EW77,V10a} for a two-well potential that is $C^1$ with a non-degenerate spinodal range show that there may be stable solutions with velocity within some intervals in the low-velocity regime;  the extent of these intervals depends on a parameter of the potential (width of the spinodal range) that is zero in our problem.  For a fully nonlinear smooth potential, analytical conclusions seem untenable at the moment; however, the numerical simulations in \cite{PK84,V10a} suggest that smoother potentials result in further stabilization of low-velocity solutions. Meanwhile, the apparent instability of the new solutions we construct below the threshold velocity suggests that the non-smooth biquadratic model cannot adequately describe stable steady propagation of a dislocation at low velocity.

The paper is organized as follows. In Section~\ref{sec:AC}  we recall  the dynamical driven  Frenkel-Kontorova model and  the equation for  traveling wave solutions that describes   steady motion of a moving dislocation, or kink, under stress.  We consider the case of a piecewise quadratic potential and recall the main properties of the  explicit Atkinson-Cabrera solution. The loss of admissibility of this solution in the low-velocity regime motivates us to relax the strict monotonicity assumption made in \cite{AC65}. We  derive in Section~\ref{sec:new_solns} conditions for a new type of kink solutions. These conditions include a linear integral equation of the first kind, whose kernel is related to  the formal Atkinson-Cabrera solution, regardless of the admissibility of the latter. Solving the integral equation yields a shape function that can be used to obtain the new kink solution. The support $[-z,z]$ of the shape function, which depends on the velocity of the moving kink, is the interval where the new kink solutions take the spinodal value. In Section~\ref{sec:small_z} we use piecewise linear and piecewise quadratic approximations of the kernel in the integral equation to obtain approximations of the shape function under the assumption that $z$ is small. It turns out  that the shape function is a distribution that involves two delta functions concentrated at $\pm z$. In contrast, the shape functions obtained in \cite{V10a} are bounded.  Our results show that the new solutions bifurcate from the Atkinson-Cabrera solutions at the velocity $V_0$ at which the latter become inadmissible. In Section~\ref{sec:numer} we describe the numerical procedure we use to obtain solutions in the case when $z$ is not necessarily small. We use this procedure to generate new solutions in the low-velocity regime $0<V<V_0$ and discuss their properties. We verify the numerical procedure by comparing it with the  analytical results of Section~\ref{sec:small_z}. In Section~\ref{sec:stab} we investigate the stability of the traveling wave solutions using numerical simulations of the discrete chain dynamics. The results suggest  instability of the new solutions.  Section \ref{sec:Visc} adds viscosity to the model. We find that this addition does not seem to stabilize slow new-type traveling waves in general.

\section{Atkinson-Cabrera traveling wave solutions}
\label{sec:AC}

The dynamics of the driven Frenkel-Kontorova chain are described by the following equation, expressed in dimensionless quantities:
\beq
\ddot{u}_n=u_{n+1}-2u_n+u_{n-1}+\mu(\sigma-\Phi'(u_n)).
\label{eq:dyn}
\eeq
Here $\sigma$ is the constant applied stress, $u_n(t)$ is the displacement of the $n$th mass at time $t$, $\mu$ is a ratio of  stiffness of the nonlinear interaction with the substrate to nearest neighbor interactions and $\Phi$ is the multiple well substrate potential. The equation describes
a Hamiltonian system with conserved energy.
To model a steadily moving dislocation, we seek solutions of \eqref{eq:dyn} in the form of a traveling wave with (constant)  velocity $V>0$:
\beq
u_n=u(\xi), \quad \xi=n-Vt.
\label{eq:TW_ansatz}
\eeq
Substituting this ansatz in \eqref{eq:dyn}, we obtain the advance-delay differential equation
\beq
V^2u''= u(\xi+1)-2u(\xi)+u(\xi-1)+\mu(\sigma-\Phi'(u(\xi))).
\label{eq:TW}
\eeq
We are interested in solutions of \eqref{eq:TW} that are of  \emph{kink type}. These satisfy the following conditions at infinity:
\beq
\langle u(\xi) \rangle \rightarrow u_\pm \quad \text{as $\xi \rightarrow \pm \infty$},
\label{eq:BCs}
\eeq
where $u_{\pm}$ are stable constant (uniform) equilibrium solutions of \eqref{eq:TW} located in two \emph{different} wells:
\[
\Phi'(u_\pm)=\sigma, \quad u_{-}>u_+, \quad \Phi''(u_\pm)>0.
\]
The angular brackets in \eqref{eq:BCs} denote the average value of the displacement because we expect this Hamiltonian discrete system
to develop oscillations. For this reason, such solutions are sometimes called \emph{generalized kinks} or heteroclinic connections of center manifolds. The average is taken over a sufficiently  large interval.

Instead of the usual periodic potential, we choose a potential with only two wells; this is appropriate for describing \emph{twinning dislocations}.
As first shown in \cite{AC65}, an explicit solution of \eqref{eq:TW}, \eqref{eq:BCs} can be obtained if one assumes
that the substrate potential is piecewise quadratic:
\beq
\Phi(u)=\begin{cases}
                          \frac{1}{2}(u+1)^2, & u \leq 0\\
                          \frac{1}{2}(u-1)^2, & u \geq 0.
           \end{cases}
\label{eq:Phi}
\eeq
Note that the derivative of this potential is discontinuous at $u=0$:
\beq
\Phi'(u)=u-2\theta(u)+1,
\label{eq:Phi_prime}
\eeq
where $\theta(u)$ is the unit step function: $\theta(u)=1$ for $u>0$, $\theta(u)=0$ for $u<0$. Observe also that in this case one has $u_\pm=\sigma \mp 1$ in
\eqref{eq:BCs}.

Suppose that the displacement switches from the second to the first well at $\xi=0$, so that
\beq
u(\xi)>0 \quad \text{for $\xi<0$,} \quad u(\xi)<0 \quad \text{for $\xi>0$}
\label{eq:AC_constraints}
\eeq
and
\beq
u(0)=0.
\label{eq:AC_switch}
\eeq
Then one may replace $\theta(u(\xi))$ by $\theta(-\xi)$ in  \eqref{eq:Phi_prime}.    Then  \eqref{eq:TW} becomes a linear equation:
\beq
V^2u''= u(\xi+1)-(2+\mu)u(\xi)+u(\xi-1)+\mu(\sigma-1+2\theta(-\xi)),
\label{eq:TW_AC}
\eeq
which can be solved using Fourier transforms. It follows from \eqref{eq:TW_AC} that the solution is of class $C^1$ at $\xi=0$ and $C^2$ elsewhere.

It should be emphasized that solutions of
\eqref{eq:TW_AC} satisfy the original nonlinear equation \eqref{eq:TW} if and only if the \emph{admissibility condition} \eqref{eq:AC_constraints}  holds. Otherwise, if a solution of \eqref{eq:TW_AC} violates \eqref{eq:AC_constraints} it will be labeled as \emph{inadmissible}.

The solution of \eqref{eq:TW_AC} constructed by Atkinson and Cabrera  \cite{AC65} (see also \cite{CB03a,KT03} for more details) is as follows:
\beq
u=U(\xi) \equiv \sigma-1-\dfrac{\mu}{i\pi}\int\limits_\Gamma\dfrac{e^{ik\xi}}{kL(k,V)}dk,
\label{eq:U_int}
\eeq
where the contour $\Gamma$ coincides with the real axis from $-\infty$ to $\infty$, except near the singular points on the real axis. One such singularity is at $k=0$, and the contour passes it from above. The other singularities are associated with the real roots of the equation $L(k,V)=0$, where
\beq
L(k,V)=\mu+4\sin^2\dfrac{k}{2}-V^2k^2.
\label{eq:L}
\eeq
The dispersion relation is the equation $L(k,\omega/k)=0$ that implicitly defines $\omega=\omega(k)$.
These roots correspond to lattice waves (phonons) emitted by the moving dislocation, and must be resolved in a way that ensures that the corresponding modes propagate ahead or behind the dislocation in accordance with the radiation condition \cite{AC65,KT03}. This condition, also known as the causality principle \cite{Slepyan02},
selects solutions such that in a frame moving with the dislocation,  lattice waves can only be emitted by the moving front and must carry energy away from it (thus causing radiative damping). Thus phonon modes whose group velocity $V_g$ is less than the velocity $V$ of the front must be placed behind it, while the modes with group velocity above $V$ propagate ahead; see \cite{KT03} for details. Setting $V=\omega(k)/k$ in \eqref{eq:L}, one shows after a little algebra that the group velocity $V_g=\omega'(k)=V+L_k(k,V)/(2Vk)$, where $L_k \equiv \partial L/\partial k$. Values of $V$ such that
 \beq L(k,V)=0\quad\hbox{and }\;L_k(k,V)=0\label{eq:reson}\eeq
for some real $k$ are the \emph{resonance velocities}. Assuming that $V>0$ is \emph{non-resonant}, i.e. $L_k(k,V) \neq 0$ at real $k$ such that $L(k,V)=0$, the sign of $V_g-V$ is that of $kL(k,V)$, and we can define
the corresponding sets of nonzero real roots as follows:
\beq
N^{\pm}(V)=\{k:L(k,V)=0,\;\mbox{\rm Im}(k)=0,\;k L_k(k,V) \gtrless 0\}.
\label{eq:Npm}
\eeq
The strict inequalities in the above definition indicate that roots at resonant velocities are excluded.
The contour $\Gamma$ goes above the roots from $N^{-}(V)$ (so that the corresponding lattice waves appear at $\xi<0$) and below the roots from $N^{+}(V)$ ($\xi>0$) in accordance with the radiation condition just mentioned. Closing the contour of integration in the upper half of the complex plane at infinity  for $\xi>0$ and lower half-plane for $\xi<0$ and applying residue theorem, one obtains (for non-resonant $V$)
\beq
U(\xi)=\begin{cases} \sigma-1-2\mu\underset{k \in M^+(V)}\sum\dfrac{e^{ik\xi}}{kL_k(k,V)}, & \xi>0\\
                       \sigma+1+2\mu\underset{k \in M^-(V)}\sum\dfrac{e^{ik\xi}}{kL_k(k,V)}, & \xi<0.
\end{cases}
\label{eq:AC_soln}
\eeq
Here
\beq
M^{\pm}(V)=\{k: L(k,V)=0,\;\mbox{\rm Im}k \gtrless 0\} \cup N^{\pm}(V)
\label{eq:Mpm}
\eeq
are the sets of all (simple) nonzero roots of the dispersion relation $L(k,V)=0$ in the corresponding half-planes.

The nonlinearity of the problem is contained in  conditions \eqref{eq:AC_constraints} and \eqref{eq:AC_switch}. Applying \eqref{eq:AC_switch} and using the continuity of $U(\xi)$ at $\xi=0$, one
obtains the \emph{kinetic relation} between the applied stress and the dislocation velocity:
\beq
\sigma=\Sigma(V) \equiv \mu\underset{k \in N(V)}\sum\dfrac{1}{|kL_k(k,V)|},
\label{eq:AC_sigma}
\eeq
where the sum is over the set of all real roots, $N(V)=N^+(V) \cup N^-(V)$; the contribution of other roots cancels out due to the symmetry of the roots. Details may be found in  \cite{KT03}. Thus the applied stress is determined entirely by the real roots. As shown in \cite{KT03}, one can derive \eqref{eq:AC_sigma} by accounting for the energy fluxes carried by the phonon waves.

 For a given non-resonant $V>0$ the solution is thus given by \eqref{eq:AC_soln}, \eqref{eq:AC_sigma}, \emph{provided} that the admissibility inequalities \eqref{eq:AC_constraints}
are satisfied.

\begin{figure}
\centerline{\psfig{figure=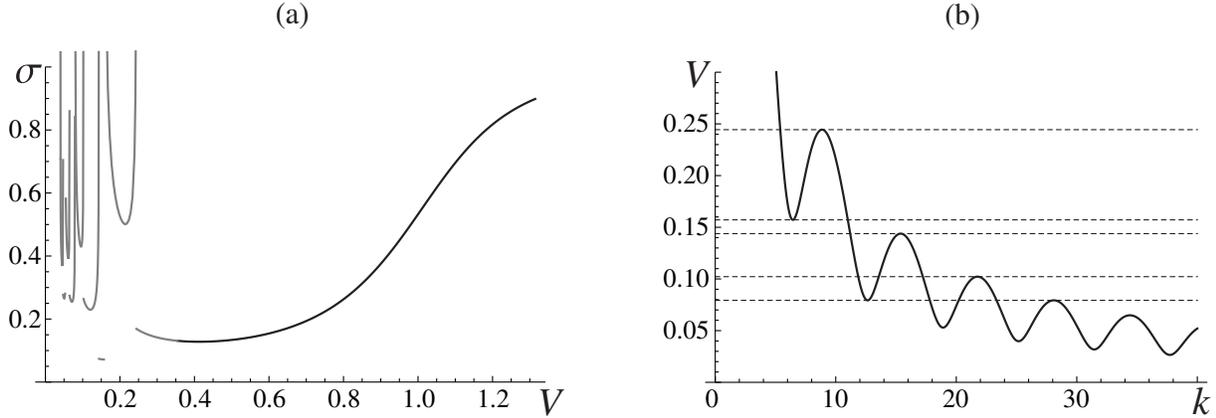,width=\textwidth}}
\caption{(a) Kinetic relation resulting from the formally obtained Atkinson-Cabrera solution.
Only the first twelve segments are shown. The gray curves correspond to inadmissible traveling waves.
(b) Solutions of $L(k,V)=0$ for positive real $k$. The dashed lines indicate the first five resonance velocities at which $L_k(k,V)=0$. Here $\mu=1$.}
\label{fig:ACkinetics}
\end{figure}

\begin{figure}
\centerline{\psfig{figure=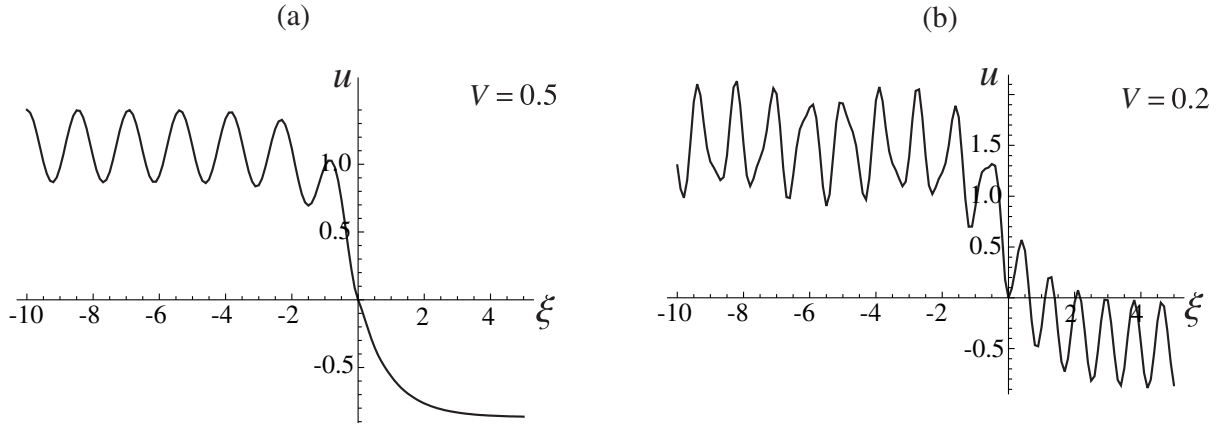,width=\textwidth}}
\caption{Displacement profiles formally computed from \eqref{eq:U_int}, \eqref{eq:AC_sigma} at (a) $V=0.5$ and (b) $V=0.2$.
Solution in (a) satisfies the constraints \eqref{eq:AC_constraints} but the one in (b) does not. Here $\mu=1$.}
\label{fig:displ_zero_z}
\end{figure}

Computing the real roots of \eqref{eq:L} for each $V>0$, we formally obtain the kinetic relation \eqref{eq:AC_sigma},
shown in Fig.~\ref{fig:ACkinetics}a for the case of $\mu=1$.
The relation consists of disjoint segments separated by resonance velocities where \eqref{eq:reson} holds (see Fig.~\ref{fig:ACkinetics}b).
A typical solution above the first resonance ($V=0.5$) is shown in
gray in Fig.~\ref{fig:displ_zero_z}a.
One can see that a moving dislocation emits phonon oscillations behind it, with  wave number equal to the single positive real root
of \eqref{eq:L}. As velocity decreases below the first resonance (see the displacement profile at $V=0.2$ in Fig.~\ref{fig:displ_zero_z}b),
more oscillation modes appear, and $u(\xi)$ formally obtained from \eqref{eq:AC_soln}, \eqref{eq:AC_sigma} features phonon emission on both sides.
However,  closer inspection
reveals that this solution of  \eqref{eq:TW_AC} is in fact \emph{inadmissible} and should be removed because it violates the assumption that $u(\xi)<0$
for $\xi>0$ (one of the admissibility conditions \eqref{eq:AC_constraints}). In fact, numerical calculations of the solution $U$ of   \eqref{eq:TW_AC} reveal that all segments of the kinetic relation below the first resonance  contain only inadmissible traveling
waves that change signs more than once,  and thus need to be removed, while the remaining large-velocity segment contains admissible solutions above a certain
threshold velocity \cite{KT03}. In the case of $\mu=1$, solutions of  \eqref{eq:TW_AC} are admissible for  $V\geq V_0 \approx 0.357$. This implies non-existence of traveling wave solutions
\emph{under the assumptions} \eqref{eq:AC_constraints}, \eqref{eq:AC_switch} when the velocity is below   the threshold value $V_0$. While this nonexistence issue was recognized  by various authors
\cite{AC65,EW77,KT03,CB03a}, it remains unclear whether there exist any moving kink type solutions below the threshold velocity.  This paper addresses this question.

\section{New traveling wave solutions}
\label{sec:new_solns}
To obtain valid traveling wave solutions of \eqref{eq:TW} in the low-velocity regime, we now replace the admissibility conditions \eqref{eq:AC_constraints}, \eqref{eq:AC_switch}
by a more general assumption
\beq
\begin{split}
&u(\xi)>0, \quad \xi<-z\\
&u(\xi)=0, \quad |\xi| \leq z\\
&u(\xi)<0, \quad \xi>z,
\end{split}
\label{eq:new_constraints}
\eeq
where $z > 0$ is to be determined. In other words, we assume that instead of switching wells instantaneously at $t=n/V$ ($\xi=0$),
the $n$th particle may stay at the spinodal value $u=0$ between the wells over the time period $[(n-z)/V,(n+z)/V]$ ($|\xi| \leq z$).
For $z=0$ we recover  Atkinson-Cabrera solutions.

Adopting a method used in \cite{FCC77,V10a}, we  observe that $\Phi'(u(\xi))$ can now be written as
\beq
\Phi'(u(\xi))=u(\xi)+1-2\int_{-z}^{z} h(s)\theta(s-\xi)ds,
\label{eq:h_1}
\eeq
where we have introduced an unknown shape function $h(s)$, which vanishes outside the interval $[-z,z]$ and is normalized so that
\beq
\int_{-z}^z h(s)ds=1.
\label{eq:h_2}
\eeq
Thus we obtain
\beq
V^2u''+(\mu+2)u(\xi)-u(\xi+1)-u(\xi-1)=\mu\biggl[\sigma-1+2\int_{-z}^z h(s)\theta(s-\xi)ds \biggr].
\label{eq:TW2}
\eeq
For consistency, we must require that in addition to \eqref{eq:TW2} and \eqref{eq:BCs}, the solution satisfies
the  \emph{generalized admissibility conditions} \eqref{eq:new_constraints}. This means that the shape function $h(\xi)$ and $z$ in \eqref{eq:h_1} must be such that these conditions are satisfied. Consistent $h(\xi)$ and $z$ then yield solutions $u(\xi)$ of \eqref{eq:TW2} that are $C^1$ at $\xi=\pm z$ and $C^2$ elsewhere.  The smoothness of solutions is discussed further in the last paragraph of Section \ref{sec:small_z}.

Taking the  Fourier transform of \eqref{eq:TW2} and using the convolution theorem, one can show (see \cite{FCC77,V10b} for details)
that
\beq
u'(\xi)=-\int_{-z}^{z} h(s)q(\xi-s)ds,
\label{eq:dudxi}
\eeq
where the kernel is the negative derivative of the solution \eqref{eq:AC_soln} of \eqref{eq:TW_AC}  with $z=0$:
\beq
q(\xi)=-U'(\xi)=\begin{cases} 2\mu i\underset{k \in M^+(V)}\sum\dfrac{e^{ik\xi}}{L_k(k,V)}, & \xi>0\\
                    -2\mu i\underset{k \in M^-(V)}\sum\dfrac{e^{ik\xi}}{L_k(k,V)}, & \xi<0,
\end{cases}
\label{eq:kernel}
\eeq
or, in the integral form (recall \eqref{eq:U_int}),
\beq
q(\xi)=\dfrac{\mu}{\pi}\int\limits_\Gamma\dfrac{e^{ik\xi}}{L(k,V)}dk.
\label{eq:q_int}
\eeq
At the same time, our assumption that $u(\xi) \equiv 0$ at $|\xi| \leq z$ implies that $u'(\xi) \equiv 0$ in the interval $(-z,z)$.
Together with \eqref{eq:dudxi}, this yields the integral equation
\beq
\int_{-z}^z h(s)q(\xi-s)ds=0, \quad |\xi|<z.
\label{eq:int_eqn}
\eeq
Thus the shape function $h(\xi)$ is an eigenfunction of the integral operator in the left hand side
of \eqref{eq:int_eqn} (with  kernel given by \eqref{eq:kernel})
associated with the zero eigenvalue. As described in more detail  in Section \ref{sec:small_z}, the integral operator has a zero eigenvalue provided that $z$ takes on special values. Note that \eqref{eq:int_eqn} is a Fredholm integral equation of the first kind.
We remark that if instead of \eqref{eq:Phi_prime} one considers the trilinear (continuous) $\Phi'(u)$ studied  in \cite{V10a},
\beq\label{eq:trilinear}
\Phi'(u)=\begin{cases}
                          u+1, & u<-\gamma/2\\
                          (1-2/\gamma)u, & |u| \leq \gamma/2\\
                          u-1, & u>\gamma/2,
           \end{cases}
\eeq
the same procedure yields a \emph{Fredholm integral equation of the second kind}
\beq\label{eq:2nd}\int_{-z}^z h(s)q(\xi-s)ds=\gamma h(\xi),\quad   |\xi|<z,
\eeq
with  the right hand side of \eqref{eq:int_eqn} replaced by $\gamma h(\xi)$,
where  the parameter $\gamma>0$  is the width of the spinodal region connecting the two wells \cite{V10a} (see also related problems with different kernels in \cite{FCC77, V10b}).
In the limit as the spinodal width parameter $\gamma \rightarrow 0+$,   the potential whose derivative is \eqref{eq:trilinear} approaches the Atkinson-Cabrera potential satisfying \eqref{eq:Phi_prime}, while  \eqref{eq:2nd} reduces to  \eqref{eq:int_eqn}.

The problem thus reduces to solving the integral equation \eqref{eq:int_eqn} for $z$ and $h(\xi)$.
Once $h(\xi)$ and $z$ are known, the convolution theorem yields the traveling wave solution:
\beq
u(\xi)=\sigma-\Sigma(V)+\int_{-z}^z h(s)U(\xi-s)ds,
\label{eq:soln}
\eeq
where the applied stress
\beq
\sigma=\Sigma(V)-\dfrac{1}{2}\int_{-z}^z h(s)(U(z-s)+U(-z-s))ds
\label{eq:sigma}
\eeq
is found by applying $u(z)=u(-z)=0$.
If $z=0$, the shape function reduces to the Dirac delta function $h(s)=\delta(s)$, and \eqref{eq:soln}, \eqref{eq:sigma} reduce to
\eqref{eq:AC_soln}, \eqref{eq:AC_sigma}, respectively.

We will present numerical evidence that this procedure yields  valid traveling waves (solutions of  \eqref{eq:TW}) for values of $V$ where the Atkinson-Cabrera solution $U$ of  \eqref{eq:TW_AC} is inadmissible. This is due to the fact that for such values of  $V$,  $u(\xi)$ given by \eqref{eq:soln} conforms to the generalized admissibility conditions \eqref{eq:new_constraints}.

\section{Kernel approximations for small $z$ and bifurcation}
\label{sec:small_z}

\begin{figure}
\centerline{\psfig{figure=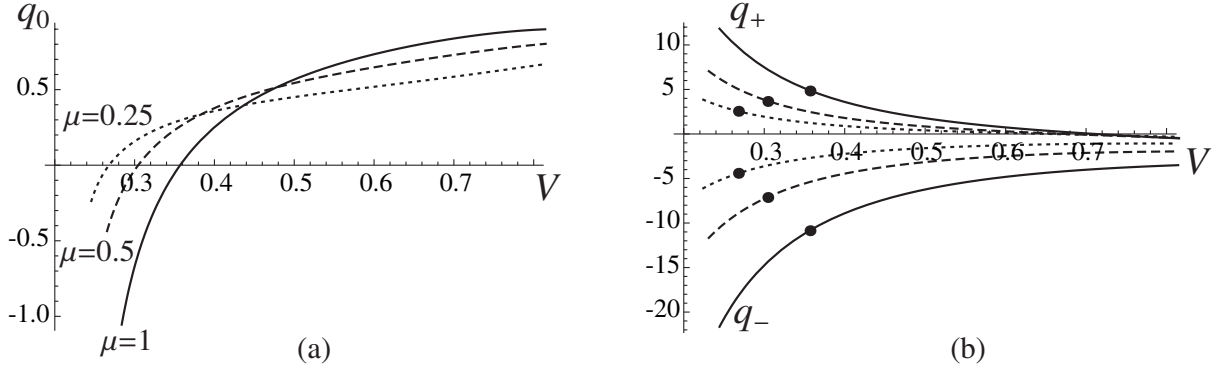,width=\textwidth}}
\caption{(a) $q_0=q(0)$ and (b) $q_{\pm}=q'(0\pm)$ as functions of $V$ for velocities above the first resonance and different values of $\mu$. The circles in (b) mark the values of $q_{\pm}$ at $V_0$ where $q_0$ changes sign.}
\label{fig:q0qplmin}
\end{figure}

Using an approximation of the integral operator in \eqref{eq:int_eqn}, we provide analytical evidence that the new type of traveling waves with $z>0$ bifurcate from the Atkinson-Cabrera solution precisely at the threshold velocity $V_0$, below which the latter becomes inadmissible. To begin with, note that in view of \eqref{eq:kernel}  we have $q(\xi)=-U'(\xi)$, where $U(\xi)$ is the Atkinson-Cabrera solution, which, as we recall, is of class $C^1$ at $\xi=0$ and $C^2$ elsewhere.  On the other hand,  the admissibility conditions \eqref{eq:AC_constraints} for $U$ clearly imply that  $U'(0)\le 0$.  Hence, \emph{a sufficient condition for loss of admissibility is $q(0)<0$}. Computing $q(0)$ as a function of $V$ numerically from \eqref{eq:kernel}, we see in Fig.~\ref{fig:q0qplmin}a  that it is monotone increasing and changes signs at a single velocity $V_0$. Consequently,  \emph{the Atkinson-Cabrera solution is inadmissible  at velocities below $V_0$}.  This argument does not in fact prove  that it is admissible above $V_0$, but  the numerical results for $U(\xi)$ at various velocities confirm that   this is the case; see also \cite{KT03}. The conclusion is:
\emph{The velocity $V_0$ at which $q(0)$ changes signs is the threshold velocity below which the Atkinson-Cabrera solution becomes inadmissible.}

As observed in the previous section, setting $z=0$ in \eqref{eq:soln}, \eqref{eq:sigma}  corresponds to the Atkinson-Cabrera solutions. Thinking of the velocity $V$ as a bifurcation parameter, it is likely that at some critical value of $V$, the new solutions with $z>0$ bifurcate off  the Atkinson-Cabrera solution. The conditions near bifurcation  are thus likely to  involve small but positive values of $z$. This suggests that in order to study the problem analytically,  we may replace the kernel $q(\xi)$ in \eqref{eq:int_eqn} by a suitable approximation near $\xi=0$, as in \cite{FCC77}.  Since $q(\xi)=-U'(\xi)$ and $U(\xi)$ satisfies \eqref{eq:TW_AC}, it follows that $q(\xi)$ is continuous, while $q'(\xi)$ has a jump discontinuity at $\xi=0$. Note also that \eqref{eq:kernel} and the continuity at $\xi=0$ imply that $q(\xi) =q_{>}(\xi)$ for $\xi\ge 0$  and  $q(\xi) =q_{<}(\xi)$ for $\xi\le 0$, where the functions $q_{>}$ and $q_{<}$, defined for $\xi\ge 0$ and $\xi\le 0$, respectively, have derivatives of at least third order in their respective domains, including $\xi=0$. Let
$$
q_0=q(0)=q_{>}(0)=q_{<}(0),
\quad q_{\pm}=q'(0\pm)=q_{\gtrless}'(0),
$$
and observe that \eqref{eq:TW_AC} implies that
\beq
\label{eq:qprime_jump}
q_{+}-q_{-}=2\mu/V^2.
\eeq
A first-order Taylor expansion of $q_{\gtrless}(\xi)$ at $\xi$ near $0$ gives $q_{\gtrless}(\xi)=q_0 + q_{\pm}  \xi + O(\xi^2)$ for $ \pm \xi\ge 0$.
Using these expansions, we obtain
 \beq
 \label{eq:pwL}
\hat{q}(\xi)=\begin{cases}q_0 + q_+ \xi ,&\xi \geq 0\\
q_0 + q_-\xi,& \xi \leq 0,\end{cases}
\eeq
 a piecewise linear approximation of $q(\xi)$ near  $\xi=0$, with an error of $O(\xi^2)$. Consider the approximate version of \eqref{eq:int_eqn} given by
 \beq
 \int_{-z}^z h(s)\hat q(\xi-s)ds=0,   \quad |\xi|<z
 \label{eq:1stkind}
 \eeq
It turns out that the only $L^1$ solution of this  is the trivial one. To see this, note that
 $$\int_{-z}^z h(s)\hat q(\xi-s)ds=\int_{-z}^\xi [q_0+q_+(\xi-s)]h(s)ds+\int_{\xi}^z [q_0+q_-(\xi-s)]h(s)ds,$$
hence it  can be differentiated twice with respect to $\xi$, with second derivative equal to $(q_+-q_-)h(\xi)$ ($ |\xi|<z$), which therefore has to vanish.  It is known that certain equations of the first kind sometimes only possess  generalized solutions involving delta functions \cite{Sakhnovich80}.  No such behavior is encountered in \cite{V10a}.

One way to approach  the possibility of generalized solutions is to consider solutions of the corresponding equation of the second kind, which actually arises in \cite{V10a},
\beq
\label{eq:sekind} \int_{-z}^z h(s)\hat q(\xi-s)ds-\gamma h(\xi)=0,\quad   |\xi|<z,
\eeq
with $\gamma$ a positive constant, and then take the limit as $\gamma\to 0+$.  We recall that this is the small-$z$ approximation of \eqref{eq:2nd}, which arises in the context of  the potential \eqref{eq:trilinear}
 considered in  \cite{V10a}. The problem studied there reduces to the present one as $\gamma\to 0+$.   Assuming $h$ is smooth enough, differentiate \eqref{eq:sekind}  twice to find that $h$ must satisfy the ODE
$$
\gamma h''(\xi)-(q_+-q_-)h(\xi) =0, \quad  |\xi|<z.
$$
 Inserting the general solution of this,
 \beq
 \label{eq:h_a}
 h(\xi)=c_1 e^{\xi/a}+c_2 e^{-\xi/a}, \quad
a=\sqrt{\gamma/(q_+-q_-)},
\eeq
 into \eqref{eq:sekind}  and evaluating the integral, we obtain an expression linear in $\xi$, which must  vanish for $ |\xi|<z$. This  yields a homogeneous linear system for the constants $c_1$ and $c_2$, which has nontrivial solutions, provided that  the corresponding determinant vanishes. This can be put into the form of a quadratic equation in the quantity  $e^{2z/a}$, whose coefficients depend on $z$ and $\gamma$.  Solving this quadratic  gives the following equations:
 \beq
 \label{eq:EE}
e^{2z\sqrt{(q_+-q_-)/\gamma}}=\frac{B_\pm(z,\gamma) }{A(z,\gamma)},
\eeq
where
\[
\begin{split} B_\pm(z,\gamma)& =
-  a  \left( {q_-}^2+ {q_+}^2\right)  \\
& \pm \sqrt{( {q_-}- {q_+})^2 \left( { q_0 }^2+  a ^2 ( {q_-}+ {q_+})^2\right)+4  { q_0 }  {q_-} ( {q_-}- {q_+})  {q_+} z+4  {q_-}^2  {q_+}^2 z^2} \end{split}
\]
 and
$$A(z,\gamma) =
  { q_0 } ( {q_-}- {q_+})+2  {q_-}  {q_+} (z-  a ),
$$
with $a$ as in \eqref{eq:h_a}.
Note that $A$, $B_\pm$ are continuous algebraic functions, which thus remain bounded  for finite values of their arguments. This provides a relation between $z$ and $\gamma$. Suppose that $z$ tends to a positive value as  $\gamma\to 0+$. Then the left hand side of \eqref{eq:EE}, hence also the right hand side, grows unbounded. This necessitates that the denominator $A\to 0$, or passing to the limit, $A(z,0)= 0$,  which  reduces to
\beq \label{eq:zL} z ={q_0(q_+-q_-)\over 2q_+ q_-}.
\eeq
For small $\gamma>0$ one can  determine the constants $c_1$,  $c_2$ and $z$ in terms of $\gamma$ by enforcing the normalization condition \eqref{eq:h_2}.  After substitution into \eqref{eq:h_a},
the latter  can be put into the form
\[
h(\xi)={q_-\cosh\left({\xi-z\over a}\right)-q_+\cosh\left({\xi+z\over a}\right) \over  a(q_--q_+)\sinh\left({2z\over a}\right)}.
\]
It is then straightforward to show that the limit as $\gamma\to0+$ is
\beq\label{eq:h_L}
h(\xi)={ q_- \over (q_--q_+)} \delta(\xi+z)-{q_+\over (q_--q_+)}\delta(\xi-z),
\eeq
in the sense of distributions (with delta functions at $\pm z$).  We conclude that in the limit as $\gamma\to 0+$, the (regular) solution of \eqref{eq:sekind} approaches the generalized solution \eqref{eq:h_L} of  the first-kind approximate equation \eqref{eq:1stkind}.  This motivates considering generalized solutions.
One can also find  \eqref{eq:h_L} directly by substituting the ansatz $h(\xi)=\alpha_+ \delta(\xi+z)+\alpha_-\delta(\xi-z)$ (with $\alpha_\pm$ unknown constants)
 into \eqref{eq:1stkind}.
This determines $\alpha_\pm$ and yields \eqref{eq:zL} and \eqref{eq:h_L}.

The value of $z$ corresponding to the above solution of the approximate equation \eqref{eq:1stkind} is given by
 \eqref{eq:zL}, whose sign we now investigate, since the possibility $z<0$ is meaningless in view of \eqref{eq:new_constraints}.  We recall  that $q_0$ and $q_\pm$ depend on $V$, and that $q_0$ is monotone increasing with a sign change at $V_0$.  One may calculate $q_\pm$
 for $V$ above the first resonance velocity $V_1$ using \eqref{eq:TW_AC}, \eqref{eq:AC_soln} and the properties of $L(k,V)$, which yield
\[
V^2 q_{\pm}=\pm\mu-\Sigma(V)(2+\mu)-4\mu \dfrac{\cos(r)}{rL_k(r,V)},
\]
where $r>0$ is the single positive real root of $L(k,V)=0$. Together with \eqref{eq:AC_sigma} this can be used to show that for $V>V_1$ one has $q_{+}>0$ if $\Sigma(V)<\mu/(\mu+4)$, which holds near $V_0$ according to our numerical results, and $q_{-}<0$; see Fig.~\ref{fig:q0qplmin}b, where
the result is shown for different values of $\mu$.
In particular, this calculation shows that in a neighborhood of $V_0$ one has $q_{+}>0$ and $q_{-}<0$.  Together with \eqref{eq:qprime_jump} this implies that $(q_+-q_-)/(2q_+ q_-)<0$  at and near $V_0$.  As a result, the right hand side of  \eqref{eq:zL} changes signs at $V_0$ (where $q_0$ vanishes) and is positive for $V$ below $V_0$.   For $V$ above  $V_0$ we recall that the Atkinson-Cabrera solution is admissible;  see the first paragraph of this section.
This provides strong analytical  evidence that \emph{bifurcation to the  new type of traveling wave with $z>0$ occurs precisely at the threshold velocity $V_0$, at which the Atkinson-Cabrera solution becomes inadmissible}. In the next section we show this numerically as well using the full kernel.

 One may improve the kernel approximation as follows. The second derivative of $ q(\xi)$  has a vanishing jump  discontinuity at $\xi=0$, as can be shown from \eqref{eq:TW_AC}; see also \cite{FCC77}.   Let
 $$q_2=q_{>}''(0)/2=q_{<}''(0)/2.$$
 Adding the quadratic term $q_2 \xi^2$ of the Taylor expansions of $q_{\gtrless}(\xi)$ to the piecewise linear kernel approximation \eqref{eq:pwL}, we obtain a piecewise quadratic approximation of $q(\xi)$ near $\xi=0$:
 $$
 \hat q(\xi)=\begin{cases}q_0 + q_+ \xi+q_2 \xi^2 ,&\xi \geq 0\\
q_0 + q_-\xi+q_2 \xi^2,& \xi \leq 0.\end{cases}
$$
The previous procedure can then be repeated with similar results. The analogue of \eqref{eq:h_L}  now takes the form
\[
\begin{split}
&h(\xi)=\alpha_+\allowbreak \delta(\xi+z)+\alpha_-\allowbreak\delta(\xi-z)+\zeta,\quad \zeta=-{2q_2\over q_+-q_-}, \\
&\alpha_\pm={q_+-q_- +4q_2 z\over 2(q_+-q_-)} \mp {q_++q_-\over 2(q_+-q_-+4q_2 z)}
\end{split}
\]
and  $z$ is a root of the quartic equation
$$(q_+-q_-)q_0+(4q_2q_0-2q_+q_-)z+4q_2(q_+-q_{-}) z^2+{32q_2^2\over 3} z^3 +{32q_2^3\over 3(q_+-q_-)}z^4=0.$$
Note that for $q_2=0$ the last two equations reduces to  \eqref{eq:h_L} and  \eqref{eq:zL}.   It turns out that only the smallest positive root $z$ of this quartic equation gives rise to traveling waves satisfying the generalized admissibility conditions \eqref{eq:new_constraints}; this root is well approximated by that of  \eqref{eq:zL}.   Moreover, setting $z=0$ in the last equation
above gives $(q_+-q_-)q_0=0$, hence the improved approximation once again predicts the bifurcation to occur at  $V_0$ where $q_0=0$.

We close this section by remarking on the smoothness of the actual traveling wave solutions $u(\xi)$ in the presence of generalized solutions of \eqref{eq:int_eqn}.   Suppose that the solution $h(\xi)$ of the latter is of generalized type, that is,
$$    h(\xi)=\alpha_+\delta(\xi+z)+\alpha_-\delta(\xi-z)+h_0(\xi)$$
where $h_0(\xi)$  is the regular part of the solution, continuous on $[-z,z]$, and $\alpha_\pm$ are constants. Then the corresponding traveling wave solution $u(\xi)$ is found from  \eqref{eq:soln}, which now gives
$$u(\xi)=\alpha_-U(\xi+z)+\alpha_+U(\xi-z)+\int_{-z}^z h_0(s)U(\xi-s)ds+\sigma-\Sigma(V)$$
By this equation and standard results on smoothness of convolutions, $u(\xi)$ has the same regularity as $U(\xi)$. Note that  $U(\xi)$ has a second derivative discontinuity at $ \xi=0$, and recall  that $u(\xi)=0$  for $|\xi|<z$ in view of \eqref{eq:new_constraints}.  It follows that $u(\xi)$ is of class $C^2$ everywhere except at $\pm z$ where $u''$ has jump discontinuities, while $u'$ is everywhere continuous. Therefore, there is no loss of regularity compared with the original Atkinson-Cabrera solution (see the remark after \eqref{eq:TW_AC})  despite allowing generalized solutions of the (auxiliary) integral equation  \eqref{eq:int_eqn}.

 \section{Numerical results}
 \label{sec:numer}
We now consider the general case where $z$ is not necessarily small.

\begin{figure}
\centerline{\psfig{figure=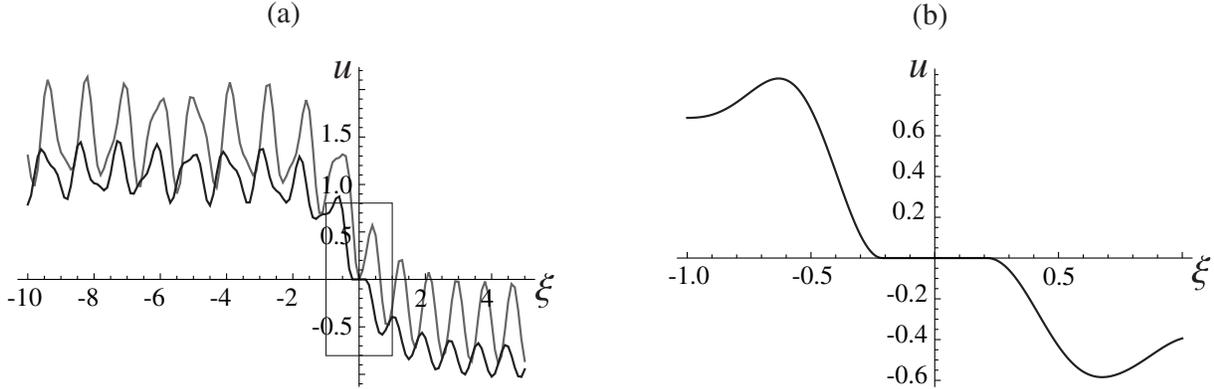,width=\textwidth}}
\caption{(a) Traveling wave solution  $u(\xi)$ at $V=0.2$, $\mu=1$ with $z=0.212$ (black curve), shown together with inadmissible $z=0$ solution  $U$ (gray curve).
(b) Zoom-in of the $z=0.212$ solution inside the rectangle in (a).}
\label{fig:displ_nonzero_z}
\end{figure}

For given non-resonant $V>0$  (such  that \eqref{eq:reson} does not hold), the solution $h(\xi)$ and $z$ of  \eqref{eq:int_eqn} was found numerically as follows. The integrals in the expressions  \eqref{eq:U_int} for the Atkinson-Cabrera solution and  \eqref{eq:q_int} for the kernel were evaluated numerically in Mathematica using Levin's method for oscillatory integrals \cite{Levin96}. The trapezoidal approximation
of the integral equation  \eqref{eq:int_eqn}  for a finite $z$ (with $100$ uniformly distributed mesh points)
then yielded a homogeneous linear system $\mathbf{Q}(z)\mathbf{h}=\mathbf{0}$. Here $\mathbf{h}$ is the vector of values of $h(s)$ at the mesh points. After solving
$\text{det}\mathbf{Q}(z)=0$ for $z$, we  obtained the corresponding solution $\mathbf{h}$ of $\mathbf{Q}(z)\mathbf{h}=\mathbf{0}$, normalized to satisfy \eqref{eq:h_2} in the sense of the trapezoidal approximation.
In general, there were more than one root of $\text{det}\mathbf{Q}(z)=0$, but at most one of these yielded admissible
solutions that satisfied the generalized admissibility conditions \eqref{eq:new_constraints} within numerical error. Once the admissible $z$ and $\mathbf{h}$ were found,
the trapezoidal approximation of the integrals in \eqref{eq:soln} and \eqref{eq:sigma} was used to compute the solution $u(\xi)$
and the applied stress $\sigma$.
Fig.~\ref{fig:displ_nonzero_z} shows the solution $u(\xi)$ obtained at $V=0.2$ and $\mu=1$ (black curve) along with the inadmissible Atkinson-Cabrera solution   $U(\xi)$
($z=0$, gray curve). In this case $z=0.212$.
One can see that the obtained solution satisfies the constraints \eqref{eq:new_constraints}
within  numerical error (which is of the order of $10^{-6}$ in this case).

\begin{figure}
\centerline{\psfig{figure=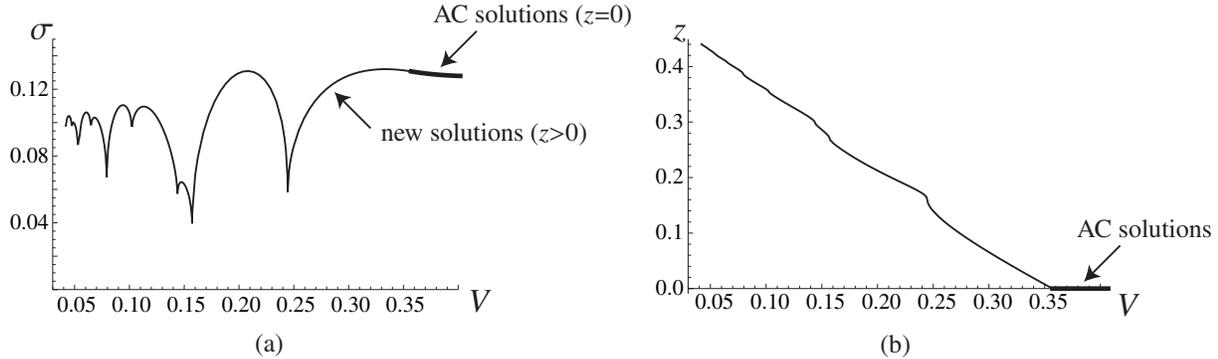,width=\textwidth}}
\caption{(a) Kinetic relation $\sigma(V)$ and (b) the corresponding $z(V)$ at $\mu=1$. The thicker segments above the threshold velocity $V_0 \approx 0.357$ indicate the parts of the curves that correspond
to Atkinson-Cabrera solutions ($z=0$).}
\label{fig:kinetics_bilinear}
\end{figure}

\begin{figure}
\centerline{\psfig{figure=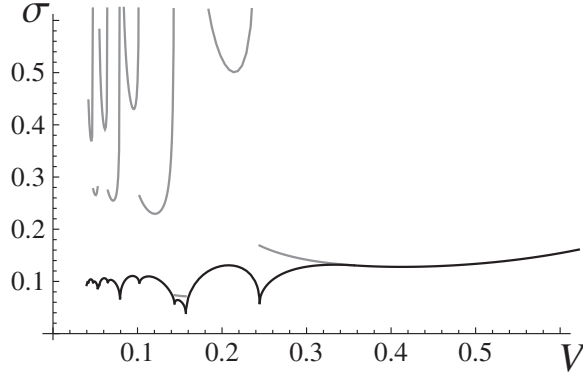,width=3.0in}}
\caption{Comparison of the kinetic relation $\sigma(V)$ (black curve) generated by the new solutions to the relation $\Sigma(V)$ (gray lines) obtained from Atkinson-Cabrera solutions. The two curves coincide above the threshold velocity $V_0 \approx 0.357$. The gray curve corresponds to \emph{inadmissible} Atkinson-Cabrera solutions below $V_0$. Here $\mu=1$.}
\label{fig:new_vs_AC}
\end{figure}

Repeating this procedure for a range of velocities, we obtain the kinetic relation $\sigma(V)$ and the corresponding relation  $z(V)$  between $z$ and $V$ shown in Fig.~\ref{fig:kinetics_bilinear}.
Note that $z(V)$ decreases as $V$ grows and becomes zero at $V=V_0 \approx 0.357$. This is the threshold velocity such that $q( 0)=0$, which equals the bifurcation velocity as predicted in Sec.~\ref{sec:small_z}.
At $V \geq V_0$, the Atkinson-Cabrera solutions are admissible, and we have $z=0$ and $\sigma=\Sigma(V)$ (thick segments in Fig.~\ref{fig:kinetics_bilinear}). Interestingly,
for $V$ below $V_0$ we obtained admissible solutions with $z>0$ even in the immediate vicinity of the \emph{resonance velocities}, where \eqref{eq:reson} holds. At these velocities, the Atkinson-Cabrera  kinetic relation
$\Sigma(V)$ is inadmissible and has singularities, while  the kinetic relation $\sigma(V)$ based on the present admissible solutions has cusps; see Fig.~\ref{fig:new_vs_AC}.

\begin{figure}
\centerline{\psfig{figure=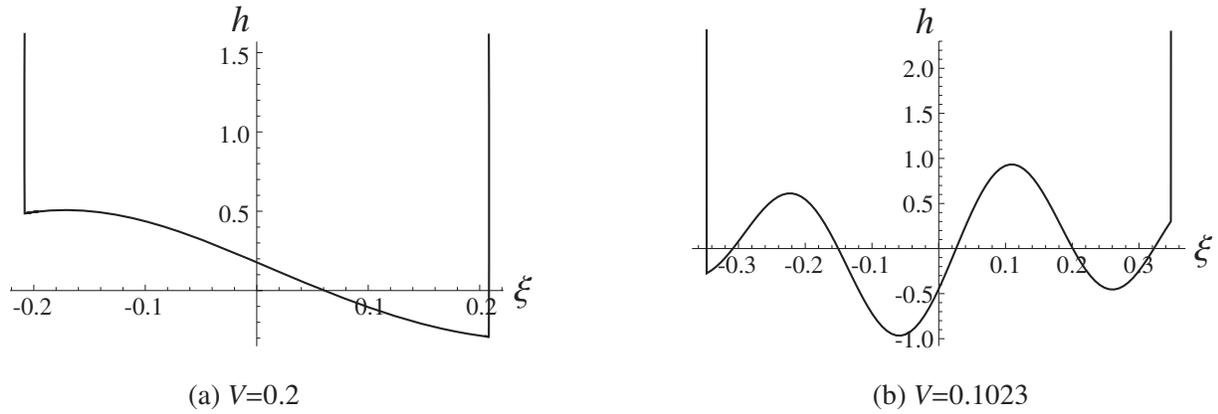,width=\textwidth}}
\caption{Shape functions $h(s)$ in the interval $[-z,z]$ at (a) $V=0.2$, $z=0.212$ and (b) $V=0.1023$, $z=0.355$. Here $\mu=1$.}
\label{fig:shape_functions}
\end{figure}

The computed shape function $h(\xi)$ corresponding to velocity $V=0.2$ ($z=0.212$) and the solution
in Fig.~\ref{fig:displ_nonzero_z}, is shown in Fig.~\ref{fig:shape_functions}a.
Note that it approximates the delta-function behavior at $\xi=\pm z$ that was predicted in Sec.~\ref{sec:small_z}. As $V$ decreases, the shape functions develop oscillations
due to the increasingly oscillatory behavior of $q(\xi)$ at smaller $V$; see, for example, Fig.~\ref{fig:shape_functions}b for $h(\xi)$ at $V=0.1023$.

We remark that one can also approximate the kernel \eqref{eq:kernel} and the Atkinson-Cabrera solution \eqref{eq:AC_soln} by replacing the infinite sums with the finite ones that include the first $N_r$ roots of \eqref{eq:L} near the origin. Using this approximation instead of the numerical evaluation of the integral, we obtained results that suggest convergence of the displacement and stress values as $N_r$ is increased. For example, at $N_r=400$, the kinetic relation is very close to the one shown in Fig.~\ref{fig:kinetics_bilinear},  the difference being of the order ranging from $10^{-5}$ for the smaller velocities to $10^{-7}$. On the other hand, the shape functions $h(\xi)$ obtained using the truncated kernel method feature additional fine-scale oscillations due to a Gibbs-type phenomenon associated with the delta-function singularities at the end points. These oscillations are probably an artifact of the truncation.

\begin{figure}
\centerline{\psfig{figure=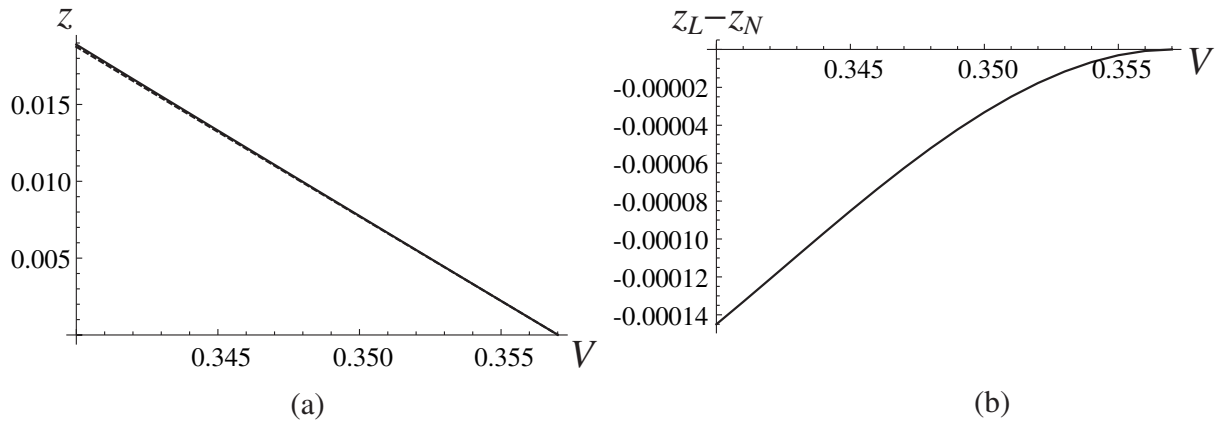,width=\textwidth}}
\caption{Comparison of the function $z_N(V)$ computed numerically using the full kernel (solid line) and $z_L(V)$ obtained from \eqref{eq:zL} using the  piecewise linear approximation \eqref{eq:pwL} of the kernel (dashed line) near the bifurcation point $V_0 \approx 0.357$. Here $\mu=1$.}
\label{fig:linear_vs_numer}
\end{figure}

Finally, we consider the bifurcation diagram for $z(V)$ near the threshold velocity $V_0$ and compare $z_N(V)$ obtained in our numerical computation to the curve $z_L(V)$ obtained from \eqref{eq:zL} using the piecewise linear approximation \eqref{eq:pwL} of the kernel (dashed line) in Fig.~\ref{fig:linear_vs_numer}.
As one can see, the two curves become closer as we approach the threshold velocity from below.

\section{Stability of traveling wave solutions}
\label{sec:stab}

To investigate the stability of the new type of traveling waves, we conducted a series
of numerical simulations, since  in general it is quite difficult to check stability of traveling waves analytically  \cite{Carpio04}. For a given applied stress $\sigma$, we used the Verlet algorithm (a symplectic scheme \cite{Verlet67}) to solve the system
\eqref{eq:dyn} of ordinary differential equation for a truncated lattice with $N$ masses, ranging from $N=600$ to $2000$, depending on the time
of the simulation. A longer chain was used if the simulation ran for a long time, in order to avoid reflection of elastic waves from the domain boundaries.
The boundary conditions were $u_0=\sigma+1$ and $u_N=\sigma-1$. Two types of initial conditions were used. The first one was Riemann-type piecewise constant initial displacement
and zero initial velocity:
\beq
(u_n(0), \dot{u}_n(0))=\begin{cases}
                          (\sigma+1, 0) & 0 \leq n<n_0\\
                          (0, 0), & n=n_0,\\
                          (\sigma-1, 0), & n_0<n \leq N
       \end{cases}
\label{eq:IC1}
\eeq
where $n_0=N/2$ for even $N$. Numerical simulations with these initial data sought to identify
stable states at a given loading that have a relatively wide basin of attraction. To capture other
possibly stable states that coexist with solutions found using \eqref{eq:IC1} but have a more narrow
basin of attraction, and to identify solutions that are likely to be unstable,
we used a second type of initial conditions, that were built from the obtained traveling wave profiles $u_n(t)=u(n-Vt)$:
\beq
(u_n(0), \dot{u}_n(0))=\begin{cases}
                          (\sigma+1, 0) & 0 \leq n<p\\
                          (u(n-n_0), -Vu'(n-n_0)) & p \leq n \leq N-p.\\
                          (\sigma-1, 0) & N-p<n \leq N
       \end{cases}
\label{eq:IC2}
\eeq
The truncated traveling wave solutions were surrounded by intervals of constant displacement of appropriately chosen size $p<n_0$
in order to ensure compatibility with the boundary conditions and avoid wave reflection from the boundaries.
In both types of simulations, after a sufficiently long time, the solution
approached an attractor corresponding to either a stationary dislocation (zero velocity) or a steadily moving front.

\begin{figure}
\centerline{\psfig{figure=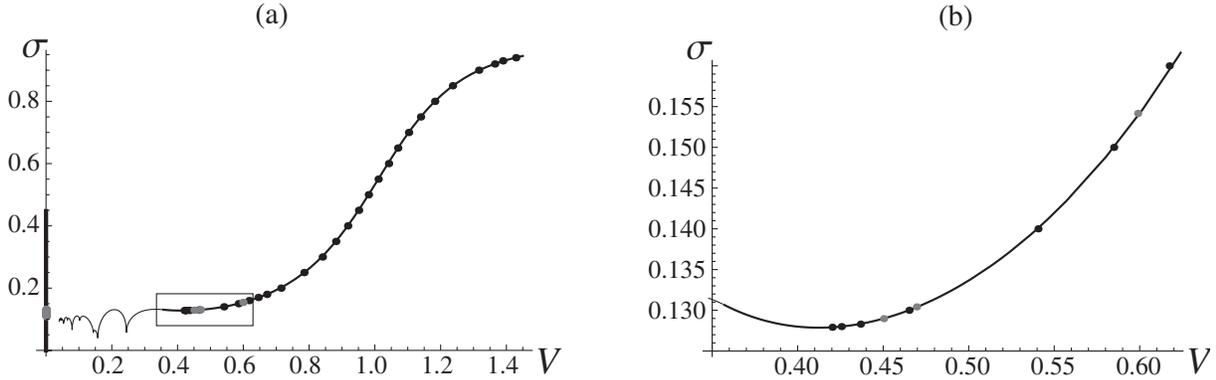,width=\textwidth}}
\caption{(a) Results of the numerical simulations at $\mu=1$ with initial data \eqref{eq:IC1} (black circles) and \eqref{eq:IC2} (gray circles),
shown together with the kinetic curve. (b) Zoom-in inside the rectangle in part (a). The thick black segment along $V=0$ axis indicates the trapping region.
The thinner portion of the kinetic curve corresponds to the solutions with $z>0$.}
\label{fig:stab_a0mu1}
\end{figure}

\begin{figure}
\centerline{\psfig{figure=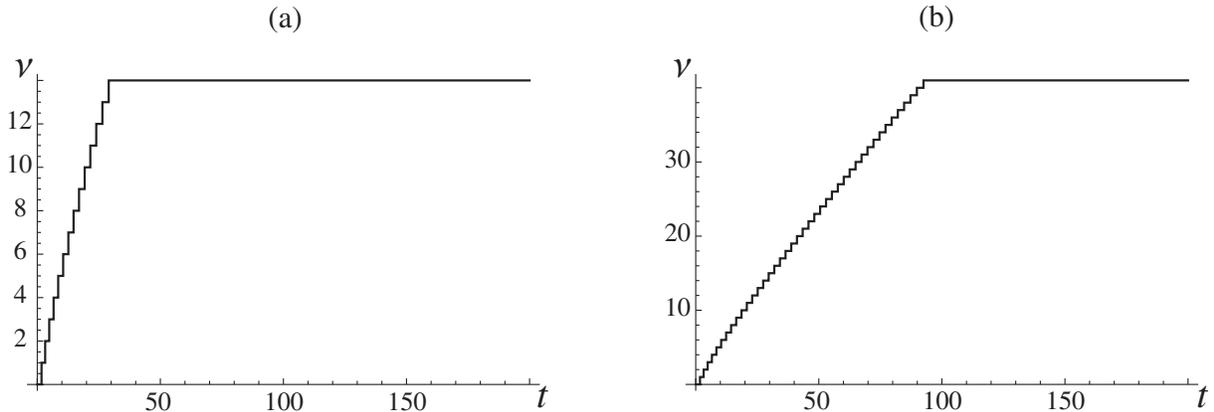,width=\textwidth}}
\caption{Position $\nu(t)$ of the dislocation at (a) $\sigma=0.125$ and (b) $\sigma=0.1275$ in the numerical simulations. Here $\mu=1$.}
\label{fig:lattice_trapping}
\end{figure}

The results are shown in Fig.~\ref{fig:stab_a0mu1}.
One can see that when the applied stress is below a certain threshold $\sigma_{\rm D}$ (here $\sigma_{\rm D}=0.128$), the long-time solution features a stationary front.
For example, in the simulations with piecewise-constant initial conditions \eqref{eq:IC1} the front propagates for some time (which increases as we
approach the $\sigma_D$ from below) and then becomes stationary. This is illustrated in Fig.~\ref{fig:lattice_trapping}, which shows
the position $\nu(t)$ of the front (defined as the value of $n$ such that $u_n(t)>0$ and $u_{n+1}(t)<0$) at $\sigma=0.125$ and $\sigma=0.1275$.
In general, stable stationary solutions exist when $\sigma$ is inside the trapping region $|\sigma| \leq \sigma_{\rm P}$, where $\sigma_{\rm P}=\sqrt{\mu/(4+\mu)}$ ($\approx0.447$ for $\mu=1$)
is the Peierls stress \cite{WS64}.  The fact that $\sigma_D<\sigma_P$ has been observed in earlier works, e.g. \cite{EW77,CB03a}. The trapping region is marked by a thick segment along $V=0$ in Fig.~\ref{fig:stab_a0mu1}.

\begin{figure}
\centerline{\psfig{figure=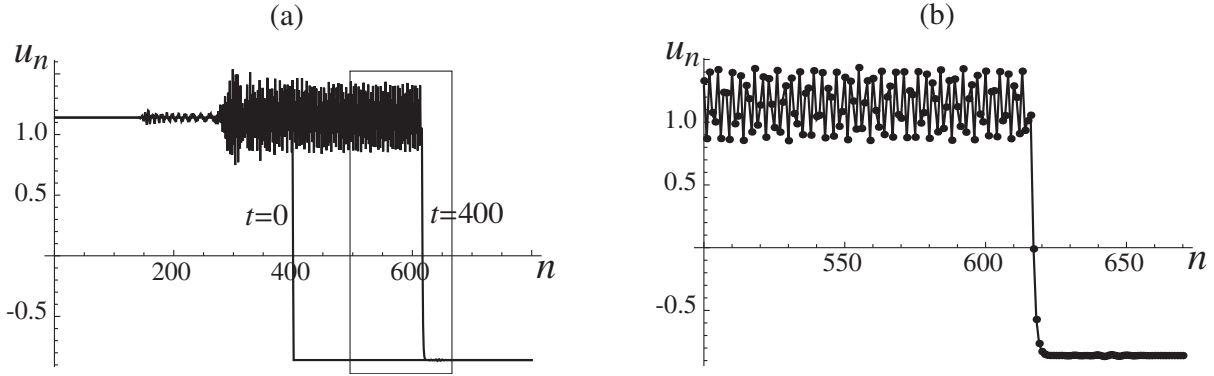,width=\textwidth}}
\caption{(a) Displacement profiles (solid lines) at $t=0$ and $t=400$ for the numerical simulation at $\mu=1$ and $\sigma=0.14$.
(b) The numerical solution (circles) at $t=600$ zoomed in around the dislocation front and
compared to the traveling wave solution (solid curve) with the same velocity, $V=0.54$.}
\label{fig:TWstability}
\end{figure}

When stress is above the threshold value ($\sigma \geq \sigma_{\rm D}$), the solution approaches a steady dislocation motion after some time. For example,
at $\sigma=0.14$, the long-time solution features the dislocation moving steadily with velocity $V=0.54$; see Fig.~\ref{fig:TWstability}.
Comparison of the numerical solution zoomed in around the front (circles) and the corresponding traveling wave solution (solid curve) in Fig.~\ref{fig:TWstability}b
shows excellent agreement, indicating stability of the traveling wave. In general, our simulations indicate that at stresses above $\sigma_{\rm D}$ all traveling waves
are stable. Note that the threshold value $\sigma_{\rm D}$ (the dynamic Peierls stress) corresponds to a local minimum of the kinetic curve, as hypothesized in \cite{AC65},
and is \emph{below} the Peierls stress $\sigma_{\rm P}$, implying that stable stationary states and stable steady motion coexist at stresses between the two values.
Similar stability results were reported in \cite{CB03a} for simulations that include a viscosity term. A proof of stability of traveling waves
with sufficiently high velocities, which for technical reasons does not extend to the whole $\sigma \geq \sigma_{\rm D}$ region, can be found in \cite{Carpio04}.

Note that inside the stability interval $\sigma \geq \sigma_{\rm D}$ suggested by the numerical simulations,
all traveling waves have $z=0$. Under  initial data \eqref{eq:IC2} based on the new-type traveling
wave solutions with $z>0$, numerical simulations always  converged to attractors with stationary fronts. Our results thus suggest that such solutions are likely to be unstable.
We remark, however, that some low-velocity solutions apparently become stable when a sufficiently wide spinodal region is included \cite{PK84,V10a}.

\section{Effect of viscosity}
\label{sec:Visc}

We now consider the effect of adding viscosity to the model on kink solutions and their stability. The rescaled governing equations become
\beq
\ddot{u}_n+\alpha\dot{u}_n=u_{n+1}-2u_n+u_{n-1}+\mu(\sigma-\Phi'(u_n)),
\label{eq:dyn_visc}
\eeq
where $\alpha>0$ is the dimensionless viscosity coefficient. Traveling wave solutions are given by \eqref{eq:AC_soln}, \eqref{eq:AC_sigma} for $z=0$ and by
\eqref{eq:soln}, \eqref{eq:sigma} for $z>0$, with \eqref{eq:L} replaced by
\beq
L(k,V)=\mu+4\sin^2\dfrac{k}{2}-V^2k^2-ik\alpha V.
\label{eq:L_visc}
\eeq
In this case there are no real roots; at small $\alpha$ the roots shift away from the real axis into upper and lower halves of the
complex plane according to the radiation condition \cite{KT03}.

\begin{figure}
\centerline{\psfig{figure=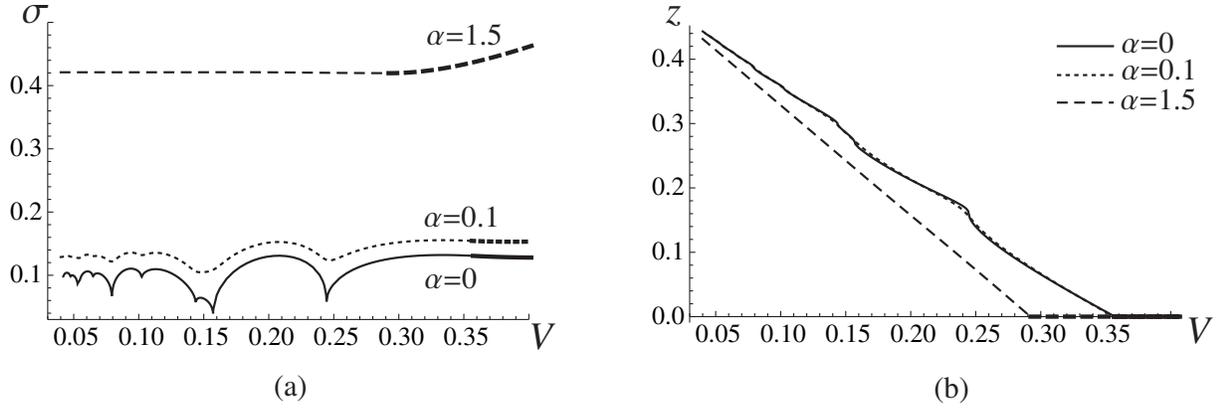,width=\textwidth}}
\caption{(a) Kinetic relation $\sigma(V)$ at $\mu=1$ and different values of the viscosity coefficient $\alpha$. The thicker portions of the kinetic curves correspond to solutions with $z=0$. (b) The corresponding $z(V)$.}
\label{fig:kinetics_with_viscosity}
\end{figure}

The effect of viscosity on the kinetic relation is shown in Fig.~\ref{fig:kinetics_with_viscosity}.
As expected, viscosity smoothens the cusps  at the resonance velocities. The amplitude of oscillations that are very pronounced in
the kinetic curves at small $\alpha$ decreases as  $\alpha$ is increased; the velocity $V_0$ at which $z=0$ solutions become admissible decreases as well. See also the corresponding
$z(V)$ graph in Fig.~\ref{fig:kinetics_with_viscosity}b. At sufficiently large $\alpha$ we have $V_0=0$, so that all $z=0$ solutions become admissible.
\begin{figure}
\centerline{\psfig{figure=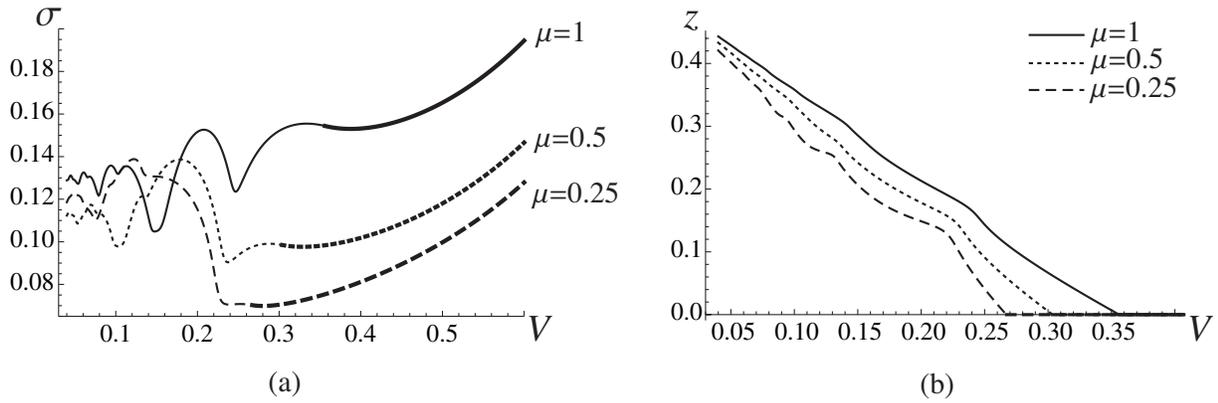,width=\textwidth}}
\caption{(a) Kinetic relation $\sigma(V)$ at $\alpha=0.1$ and different values of $\mu$. The thicker portions of the kinetic curves correspond to solutions with $z=0$. (b) The corresponding $z(V)$.}
\label{fig:effect_of_mu}
\end{figure}

The value of $V_0$ also becomes smaller as we decrease $\mu$, as shown in Fig.~\ref{fig:effect_of_mu}, where $\alpha=0.1$ (see also Fig.~\ref{fig:q0qplmin}a, where $\alpha=0$). Another effect of
the parameter $\mu$ is the different values of the resonance velocities.

\begin{figure}
\centerline{\psfig{figure=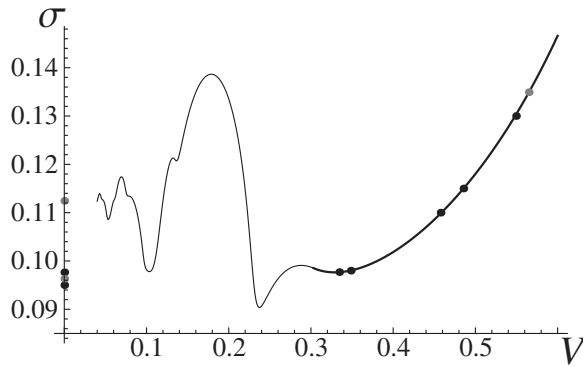,width=3.in}}
\caption{Results of the numerical simulations at $\mu=0.5$ and $\alpha=0.1$ with initial data \eqref{eq:IC1} (black circles) and \eqref{eq:IC2} (gray circles),
shown together with the kinetic curve. The thinner portion of the kinetic curve corresponds to the solutions with $z>0$.}
\label{fig:stab_a0p1mu0p5}
\end{figure}

The stability picture is not significantly affected by $\alpha$ and $\mu$. See, for example, Fig.~\ref{fig:stab_a0p1mu0p5},
where $\mu=0.5$ and $\alpha=0.1$.
As before, only sufficiently fast $z=0$ solutions, above the last local minimum of the kinetic curve, appear to be stable.

\subsection*{Acknowledgements}  This research was partially supported by  the European Union's Seventh Framework Programme (FP7-REGPOT-2009-1) under grant  no. 245749 through  the Archimedes Center for Modeling, Analysis and Computation (ACMAC) of the Department of Applied Mathematics at the University of Crete. The work of A.V. was supported by the U.S. National Science Foundation through  grant DMS-1007908. A.V. thanks ACMAC and the Department of Applied Mathematics at the University of Crete for hospitality.


\begin{thebibliography}{10}

\bibitem{AC65}
W.~Atkinson and N.~Cabrera.
\newblock Motion of a {F}renkel-{K}ontorova dislocation in a one-dimensional
  crystal.
\newblock {\em Phys. Rev. A}, 138(3):763--766, 1965.

\bibitem{Carpio04}
A.~Carpio.
\newblock Nonlinear stability of oscillatory wave fronts in chains of coupled
  oscillators.
\newblock {\em Phys. Rev. E}, 69:046601, 2004.

\bibitem{CB03a}
A.~Carpio and L.~L. Bonilla.
\newblock Oscillatory wave fronts in chains of coupled nonlinear oscillators.
\newblock {\em Phys. Rev. E}, 67:056621, 2003.

\bibitem{CF}
V.~Celli and N.~Flytzanis.
\newblock Motion of a screw dislocation in a crystal.
\newblock {\em J. Appl. Phys.}, 41(11):4443--4447, 1970.

\bibitem{EW77}
Y.~Y. Earmme and J.~H. Weiner.
\newblock Dislocation dynamics in the modified {F}renkel-{K}ontorova model.
\newblock {\em J. Appl. Phys.}, 48(8):3317--3341, 1977.

\bibitem{FCC77}
N.~Flytzanis, S.~Crowley, and V.~Celli.
\newblock High velocity dislocation motion and interatomic force law.
\newblock {\em J. Phys. Chem. Solids}, 38:539--552, 1977.

\bibitem{IS71}
S.~Ishioka.
\newblock Uniform motion of a screw dislocation in a lattice.
\newblock {\em Journal of the Physical Society of Japan}, 30(2):323--327, 1971.

\bibitem{KT03}
O.~Kresse and L.~Truskinovsky.
\newblock Mobility of lattice defects: discrete and continuum approaches.
\newblock {\em J. Mech. Phys. Solids}, 51:1305--1332, 2003.

\bibitem{Levin96}
D.~Levin.
\newblock Fast integration of rapidly oscillatory functions.
\newblock {\em Journal of Computational and Applied Mathematics}, 67:95--101,
  1996.

\bibitem{NMR11}
V.~Nosenko, G.~E. Morfill, and P.~Rosakis.
\newblock Direct experimental measurement of the speed-stress relation for
  dislocations in a plasma crystal.
\newblock {\em Phys. Rev. Lett.}, 106(15):155002--155006, 2011.

\bibitem{PK84}
M.~Peyrard and M.~D. Kruskal.
\newblock Kink dynamics in the highly discrete sine-{G}ordon system.
\newblock {\em Physica D}, 14:88--102, 1984.

\bibitem{Sakhnovich80}
L.~A. Sakhnovich.
\newblock Equations with a difference kernel on a finite interval.
\newblock {\em Russian Math Surveys}, 35(4):81--152, 1980.

\bibitem{Slepyan02}
L.~I. Slepyan.
\newblock {\em Models and phenomena in Fracture Mechanics}.
\newblock Springer-Verlag, New York, 2002.

\bibitem{V10a}
A.~Vainchtein.
\newblock Effect of nonlinearity on the steady motion of a twinning
  dislocation.
\newblock {\em Physica D}, 239:1170--1179, 2010.

\bibitem{V10b}
A.~Vainchtein.
\newblock The role of spinodal region in the kinetics of lattice phase
  transitions.
\newblock {\em J. Mech. Phys. Solids}, 58(2):227--240, 2010.

\bibitem{Verlet67}
L.~Verlet.
\newblock Computer ``experiments" on classical fluids. {I. Thermodynamical}
  properties of {Lennard-Jones} molecules.
\newblock {\em Phys. Rev.}, 159:98--103, 1967.

\bibitem{WS64}
J.~H. Weiner and W.~T. Sanders.
\newblock Peierls stress and creep in a linear chain.
\newblock {\em Phys. Rev.}, 134(4A):1007--1015, 1964.

\end{thebibliography}
\end{document}